\documentclass{elsart}

\journal{}

\usepackage{amsmath}
\usepackage{amssymb}
\usepackage{longtable}

\begin{document}

\begin{frontmatter}

 \title{A note on the zeros and local extrema of Digamma related functions}
 \author[mezo]{Istv\'an Mez\H{o}}
 \ead{istvanmezo81@gmail.com}
 \address{School of Mathematics and Statistics,\\Nanjing University of Information Science and Technology,\\No.219 Ningliu Rd, Pukou, Nanjing, Jiangsu, P. R. China}
 \thanks[mezo]{The research of Istv\'an Mez\H{o} was supported by the Scientific Research Foundation of Nanjing University of Information Science \& Technology, the Startup Foundation for Introducing Talent of NUIST. Project no.: S8113062001, and the National Natural Science Foundation for China. Grant no. 11501299.}
 
\begin{abstract}Little is known about the zeros of the Digamma function. Establishing some Weierstrassian infinite product representations for a given regularization of the Digamma function we find interesting sums of its zeros. In addition, we study the same questions for the zeros of the logarithmic derivative of the Barnes $G$ function. At the end of the paper we provide rather accurate approximations of the hyperfactorial where, rather interestingly, the Lambert function appears.
\end{abstract}

\end{frontmatter}

\section{Introduction}

The \emph{Euler gamma function} is defined by the improper integral
\[\Gamma(z)=\int_0^\infty e^{-z}t^{z-1}dt\quad(\Re(z)>0).\]
The logarithmic derivative of $\Gamma$ is the \emph{Digamma function}
\[\psi(z)=\frac{\Gamma'(z)}{\Gamma(z)}\quad(z\in\mathbb C\setminus\{0,-1,-2,\dots\}).\]
This function is analytic everywhere except the non positive integers, and it has first order poles in these points. It is known that the Digamma function has only real and simple zeros, and all of them except only one are negative.



The Barnes $G$ function \cite{SCbook} is a higher order analogue of the $\Gamma$ function satisfying the functional equation
\[G(z+1)=\Gamma(z)G(z)\]
with the normalization $G(1)=1$.

The derivative of $\log G$ does not give us a very different function, because
\[\psi_G(z)=\frac{G'(z)}{G(z)}=\frac{\log(2\pi)+1}{2}-z+(z-1)\psi(z)\quad(z\in\mathbb C\setminus\{0,-1,-2,\dots\}).\]
This function, similarly to the Digamma function, has first order poles in the non positive integers and has real zeros of multiplicity one, two of the zeros are positive.

In this paper we study a regularization of the $\psi$ and $\psi_G$ functions which gives us the possibility to find interesting closed form relations for some infinite sums of the zeros of these functions.

\section{Zeros of the Digamma function}

As we mentioned, the $\psi$ function has only real roots, and has first order poles at the non positive integers. The zeros are of multiplicity one. Canceling the poles somehow, we can get a Weierstrass product representation which is so useful many times. Especially, we can use it to find various infinite sums for the zeros. What we have found is contained in the following theorem.

\begin{thm}For all $z\in\mathbb C$ we have that
\begin{equation}
\frac{\psi(z)}{\Gamma(z)}=-e^{2\gamma z}\prod_{k=0}^\infty\left(1-\frac{z}{\alpha_k}\right)e^{\frac{z}{\alpha_k}}.\label{WPR}
\end{equation}
Here $\gamma$ is the Euler constant, and $\alpha_k$ is the $k$th zero of $\psi$ counted from the right. So, in special, $\alpha_0=1.461632\dots$ is the unique positive zero of $\psi$.
\end{thm}

\textit{Proof}. By the properties of the poles of $\psi$ and $\Gamma$ it comes that $\frac{\psi}{\Gamma}$ is an entire function and the set of zeros of this function is the same as that of $\psi$. Moreover, it is known due to Hermite \cite[p. 259]{AS} that
\begin{equation}
\alpha_k=-k+O\left(\frac{1}{\log k}\right).\label{Hermite}
\end{equation}
Hence, considering the theory of entire functions \cite{Levin} we can be sure that
\[\frac{\psi(z)}{\Gamma(z)}=ce^{az}\prod_{k=0}^\infty\left(1-\frac{z}{\alpha_k}\right)e^{\frac{z}{\alpha_k}},\]
for some constants $a$ and $c$, see \cite[p. 26]{Levin}. Since $\frac{\psi(0)}{\Gamma(0)}=-1$ we get that $c=-1$. Taking the logarithmic derivative on both sides and substituting zero again it also comes that $a=2\gamma$.\hfill\qed

Some consequences immediately follow. Taking logarithmic derivative on both sides of \eqref{WPR} we get that
\[\frac{\psi'(z)}{\psi(z)}-\psi(z)=2\gamma-z\sum_{k=0}^\infty\frac{1}{\alpha_k^2-\alpha_kz}.\]
With $z=1$ this becomes
\[\sum_{k=0}^\infty\frac{1}{\alpha_k^2-\alpha_k}=\gamma+\frac{\pi^2}{6\gamma}.\]
Another consequence of \eqref{WPR} is that
\[\frac{\psi(z)}{\Gamma(z)}\frac{\psi(-z)}{\Gamma(-z)}=\prod_{k=0}^\infty\left(1-\frac{z^2}{\alpha_k^2}\right),\]
so, comparing the coefficient of $z^2$ on both sides we have
\[\sum_{k=0}^\infty\frac{1}{\alpha_k^2}=\gamma^2+\frac{\pi^2}{2}.\]
The logarithmic derivative of the penultimate formula shows that
\[\frac{1}{z}+\pi\cot(\pi z)+\frac{\psi_1(z)}{\psi(z)}-\frac{\psi_1(-z)}{\psi(-z)}=-2z\sum_{k=0}^\infty\frac{1}{\alpha_k^2-z^2},\]
so, in special,
\[\sum_{k=0}^\infty\frac{1}{\alpha_k^2-1}=\frac{\gamma}{2}+\frac{\pi^2}{12\gamma}-1.\]

One more series identity we can deduce with small effort. Letting
\begin{equation}
F(z)=\frac{\psi(z)}{\Gamma(z)}\frac{\psi(-z)}{\Gamma(-z)},\label{F}
\end{equation}
it follows that
\[F(z)F(iz)=\prod_{k=0}^\infty\left(1-\frac{z^4}{\alpha_k^4}\right).\]
Comparing the coefficients of $z^4$ we arrive at the following identity.
\[\sum_{k=0}^\infty\frac{1}{\alpha_k^4}=\gamma^4+\frac{2\gamma^2\pi^2}{3}+\frac{\pi^4}{9}+4\gamma\zeta(3).\]

Worth to note that all of these series converge slowly.

\section{Zeros of the $\psi_G$ function}

Since the zeros of the $\psi_G$ function behave similarly to the zeros of the $\psi$ function (see the next subsection) we can repeat all the above arguments to have that
\[\frac{\psi_G(z)}{\Gamma(z)}=e^{\left(2\gamma+\frac12\log2\pi-\frac12\right)z}\prod_{k=0}^\infty\left(1-\frac{z}{\beta_k}\right)e^{\frac{z}{\beta_k}},\]
where $\beta_k$ is the $k$th zero of $\psi_G$ on the real line such that
\[\beta_0=2.55766\dots,\quad\beta_1=1.39147\dots,\quad\beta_2=-0.3662934,\]
and so on.

It follows from the Weierstrass product that
\[\sum_{k=0}^\infty\frac{1}{\beta_k^2}=\frac94+\frac{\pi^2}{2}+\gamma(1+\gamma+\log2\pi)-\frac12\log2\pi+\frac14\log^22\pi.\]
Similarly as above (see \eqref{F}), we can define the $F_G$ function to calculate the sum of the reciprocals of $\beta_k^4$. This is a bit cumbersome, but the result is
\[\sum_{k=0}^\infty\frac{1}{\beta_k^4}=2 \zeta (3) (L+1)+\gamma ^4+\frac{\pi ^4}{9}+\frac{1}{6} \pi ^2 \left(L^2+3\right)+2 \gamma ^3 (L+1)+\]
\[\frac{1}{6} \gamma  \left[24 \zeta (3)+3+4 \pi ^2 (L+1)+3L (L^2-L+7)\right]+\]
\[\frac{1}{6} \gamma ^2 \left[9 L^2+6 L+4 \pi ^2+21 \right]+\frac{1}{16}\left[L^4 -4 L^3+22 L^2-36 L+49\right],\]
where $L=\log2\pi$.

\subsection{The asymptotic behavior of the zeros of $\psi_G$}

As we already noted (see \eqref{Hermite}), Hermite knew that for the roots of the $\psi$ function behaves like $\alpha_k=-k+O\left(\frac{1}{\log k}\right)$. This can be strengthened taking into account that
\[\psi(1-x)-\psi(x)=\frac{\pi}{\tan(\pi x)},\]
substituting $x=\alpha_k$ and then taking the asymptotic approximation
\[\psi(x)=\log x-\frac{1}{2x}+O\left(\frac{1}{x^2}\right)\]
for large $x$. This argument results that
\[\alpha_k\approx-k+\frac{1}{\pi}\arctan\left(\frac{\pi}{\log k-\frac{1}{2k}+O\left(\frac{1}{k^2}\right)}\right).\]
More terms from the approximation of $\psi$ give more precise expression. As $k$ grows, the argument of $\arctan$ is close to zero. Close to the origin $\arctan(x)\approx x$, therefore Hermite's estimation follows. (We note that the idea comes from the WikiPedia article on the Digamma function but there is no reference there).

Such a trick can be applied for $\psi_G$, mutatis mutandis. The functional equation for $\psi_G$ is
\[\psi_G(1-x)+\psi_G(x)=\log2\pi-\frac{\pi x}{\tan(\pi x)}-\psi(x).\]
Substituting $x=\beta_k$ and supposing that $\psi(\beta_k)\approx 0$ (or at least converges to zero as $k\to\infty$ which is, in fact, true) we can apply the asymptotic \cite[p. 40]{SCbook}
\[\psi_G(1+x)=\frac12\log2\pi-x+x\log x+O\left(\frac1x\right)\quad(x\to\infty)\]
to have that
\[\beta_{k+2}=-k+\frac{1}{\pi}\arctan\left(\frac{\pi}{\log k-1-\frac{\log2\pi}{2k}+O\left(\frac{1}{k^2}\right)}\right).\]
In Hermitean form:
\[\beta_{k+2}=-k+\frac{1}{\log k-1-\frac{\log2\pi}{2k}+O\left(\frac{1}{k^2}\right)}.\]
(We need to shift the indices of $\beta$ by 2 to match the real location, $\beta_0$ and $\beta_1$ are positive.)

This shows that the zeros of $\psi_G$ converge a bit more slowly to the bounding left integer than $\alpha_k$ does. For example,
\[\alpha_{10}\approx-9.702672541,\quad\alpha_{100}\approx-99.80953650,\quad\alpha_{1000}\approx-999.8641415,\]
while
\[\beta_{11}\approx-9.622785495,\quad\beta_{101}\approx-99.77177415,\quad\beta_{1001}\approx-999.8444267.\]





\section{Extremal values of the hyperfactorial functions}

There are two hyperfactorial functions, both related to the Gamma and Barnes $G$ functions. One of the hyperfactorials is denoted by $K(x)$ and defined as
\[K(x)=\frac{\Gamma(x)^{x-1}}{G(x)}.\]
Named so because its values at positive integers are
\[K(n+1)=\prod_{k=1}^nk^k\quad(n\ge0).\]

Say something about the zeros and poles of this function...

With respect to the extrema of $K(x)$, we can say the following.

\begin{thm}For real $x$ the function $K(x)$ has two positive real extrema at
\[x_1=1.53769\;\;\mbox{is a local minimum},\]
\[x_0=0.290957\;\;\mbox{is a local maximum},\]
and has infinitely many extrema on the negative real axis approximately at
\[x_{-n}\sim-n+\frac{W\left(a_n(1+\log n)\right)}{1+\log n}\quad(n=1,2,3,\dots),\]
where
\[a_n=\frac{\exp\left(-\frac12(1+\log(2\pi))\right)}{\pi}\frac{\cos(\pi n)}{n^n}.\]
\end{thm}

\textit{Proof}. Taking the derivative of $K(x)$, it is easy to see that the extrema can occur at $x$ satisfying the transcendental equation
\[\log\Gamma(x)+x=\frac12+\frac12\log(2\pi).\]
We are looking for solutions on the negative real axis, so we rewrite this equation as
\begin{equation}
\log\Gamma(-x)-x=c,\label{treq}
\end{equation}
where we introduced the constant $c=\frac12+\frac12\log(2\pi)$ and we suppose that $x>0$.

Applying the reflexion formula
\[\log\Gamma(-x)+\log\Gamma(x)=\log(-\pi)-\log x-\log\sin(\pi x)\]
together with the asymptotics
\[\log\Gamma(x)=x\log x-x+\frac12\log(2\pi)+O(\log(x)),\]
we get that the exact equation \eqref{treq} transforms into the approximative equation
\[\log\csc(\pi x)-x\log x=c,\]
where $\csc=1/\sin$ and $x>0$. Taking exponential, we have that this takes the equivalent form
\begin{equation}
x^x\sin(\pi x)=e^{-c}.\label{treq2}
\end{equation}
Numerical approximations show that the zeros of this equation (and thus the extrema of $K(x)$) are very close to the negative integers, so we apply Taylor's approximation $\sin(\pi x)\approx\pi\cos(\pi n)(x+n)$, and $(x+\varepsilon)^{x+\varepsilon}\approx x^x\exp(\varepsilon(1+\log x))$ together with the assumption that $x=n+\varepsilon_n$. Hence \eqref{treq2} turns to be
\[\frac{e^{-c}}{\pi}\frac{\cos(\pi n)}{n^n}=\varepsilon\exp(\varepsilon(1+\log x))=\varepsilon(e\cdot n)^\varepsilon.\]
This equation can be solved in terms of the Lambert $W(x)$ function which gives the solutions of the equation $we^w=x$. An algebraic transformation helps to see this more directly, because if
\[\varepsilon(e\cdot n)^\varepsilon=a,\]
then
\[\varepsilon=\frac{W(a(1+\log n))}{1+\log n}.\]
From here the result follows.

\end{document}